\input amstex
\magnification=\magstep1
\input epsf
\input amssym.def
\input amssym
\pageno=1
\baselineskip 14 pt
\def \pop#1{\vskip#1 \baselineskip}

\font\gr=cmbx12

\def \rad{\operatorname {rad}}
           
\def \et{\operatorname {et}}
\def \dim{\operatorname {dim}}

\def \Spec{\operatorname {Spec}}
\def \Spf{\operatorname {Spf}}
\def \lim{\operatorname {lim}}

\def \Spm  {\operatorname {Spm  }}

\def \ord{\operatorname {ord }}
\def \ord{\operatorname {ord }}

\def \mod{\operatorname {mod }}

\def \Spf{\operatorname {Spf}}

\def \gcd{\operatorname {gcd}}

\def \SS-Deg{\operatorname {SS-Deg}}
\def \SNS-Deg{\operatorname {SNS-Deg}}
 \def \DS-Deg{\operatorname {DS-Deg}}
\def \DNS-Deg{\operatorname {DNS-Deg}}

\pop {2}
\par
\noindent                                          
\centerline {\bf \gr Wild ramification and a vanishing cycles formula}

\pop {3}
\noindent                                          
\centerline {\bf \gr Mohamed Sa\"\i di}

\pop {4}
\par
\noindent
\centerline {\bf \gr  Abstract}
\rm 
\par 
In this paper we prove an explicit formula which compares the dimensions of 
the spaces of vanishing cycles in a Galois cover of degree $p$ between formal 
germs of curves over a complete discrete valuation ring of inequal 
characteristics $(0,p)$. This formula can be easily generalised to the 
case of a Galois cover with group which is nilpotent or 
which has a normal $p$-sylow subgroup. The results of this paper 
play a key role in [Sa-1] where is studied the semi-stable reduction of 
Galois covers of degree $p$ above semi-stable curves over a complete 
discrete valuation ring of inequal characteristics $(0,p)$, as well as 
the Galois action on these covers.

\pop {2}
\par
\noindent
{\bf \gr 0. Introduction.}\rm \ Let $R$ be a complete discrete valuation ring
of inequal characteristics, with uniformiser $\pi$, fraction field $K$, 
and algebraically closed residue field $k$ of characteristic $p$. 
In this paper we investigate Galois covers of degree $p$ 
between formal germs of $R$-curves at closed points. Our main 
result is the following formula which compares the dimensions
of the spaces of vanishing cycles at the corresponding closed points. More
precisely we have the following:

\pop {.5}
\par
\noindent
{\bf \gr Theorem (2.4).}\rm \ {\sl Let $\Cal X:=\Spf {\hat {\Cal O}_x}$ be 
the formal germ of an $R$-curve at a closed point $x$, with $\Cal X_k$ 
reduced. Let $f:\Cal Y\to \Cal X$ be a 
Galois cover of group $\Bbb Z/p\Bbb Z$ with $\Cal Y$ normal and local. 
Assume that the special fibre  $\Cal Y_k:=\Cal Y\times _Rk$ of $\Cal Y$ 
is reduced. Let $\{\wp_i\}_{i\in I}$ be the minimal prime ideals of 
$\hat {{\Cal O}}_x$ which contain $\pi$ and which correspond to the  
branches $(\eta _i)_{i\in I}$ of the special fibre 
$\Cal X_k:=\Cal X\times _Rk$ of $\Cal X$ at $x$, 
and let $\Cal X_{i}:=\Spf \hat \Cal O_{\wp _i}$ be the formal completion of 
the localisation of $\Cal X$ at $\wp _i$. For each $i\in I$, the above cover
$f$ induces a torsor $f_i:\Cal Y_i\to \Cal X_{i}$ under 
a finite and flat $R$-group scheme $G_i$ of rank $p$, above the boundary 
$\Cal X_i$. For 
each $i\in I$, let $(G_{k,i},m_i,h_i)$ be the
reduction type of $f_i$ as defined in [Sa] 3.2. Let $y$ be the closed point 
of $\Cal Y$. Then one has the following {\bf ``local Riemman-Hurwitz 
formula''}:
$$\ \ \ 2g_y-2=p(2g_x-2)+d_{\eta}-d_s$$
Where $g_y$ (resp. $g_x$) denotes the genus of the singularity at $y$ 
(resp. $x$), $d_{\eta}$ is the degree of the divisor of ramification in 
the morphism
$f_{\eta}:\Cal Y_{\eta}\to \Cal X_{\eta}$ induced by $f$ on the generic 
fibre, and $d_s:=
\sum _{i\in I^{\rad}}
(m_i-1)(p-1)+\sum _{i\in I^{\et}}(m_i-1)(p-1)$, where $I^{\rad}$ is the 
subset of $I$ consisting of those $i$ for which $G_{k,i}$ is radicial, and 
$I^{\et}$ is the subset of $I$ consisting of those $i$ for which 
$G_{k,i}$ is \'etale and $m_i\neq 0$. Hier $G_{i,k}$ denotes the special fibre of the group scheme $G_i$}. 

\pop {.5}
\par
In particular the genus $g_y$ of $y$ depends only on the genus $g_x$ of $x$, 
the ramification datas on the generic fibre in the above morphism 
$f:\Cal Y\to \Cal X$, and its 
degeneration type on the boundaries of the formal fibre $\Cal X$. 
The above formula can be easily extended to the case of a Galois 
cover with group $G$ which is nilpotent, or a group which has a normal 
$p$-Sylow subgroup. Our method to prove such a formula is to construct,
using formal patching techniques \`a la Harbater, a 
compactification $\Tilde f:Y\to X$ of the above 
cover $f:\Cal Y\to \Cal X$ (cf. 2.3.2). The formula follows then by 
comparing the genus of the special and the generic fibres of $Y$ in this 
compactification.

\pop {.5}
\par
As an application, and using the above formula one can obtain interesting 
results in the case where $\Cal X$ is the formal germ of a semi-stable 
$R$-curve (cf. 3.1 and 3.2), in particular one can predict in this case 
if $\Cal Y$ is semi-stable or not. We give several examples 
which illustrate this 
situation namely the case of Galois covers of degree $p$ between 
formal germs of semi-stable $R$-curves (cf. 3.1.3, 3.1.4 and 3.2.4). 
In particular one can classify \'etale Galois covers of degree $p$ between 
annuli (cf. 3.2.5). The above results play a key role in [Sa-1] 
in order to exhibit and realise the ``degeneration datas'' 
associated to Galois covers of degree $p$ above a proper semi-stable 
$R$-curve.

\pop {1}
\par
\noindent
{\bf \gr I. Formal and rigid patching.}\rm \ 

\pop {.5}
\par
In what follows we explain 
the procedure which allows to construct 
(Galois)-covers of curves in the setting of
formal or rigid geometry by glueing together covers of 
formal affine or affinoid rigid curves with covers of formal fibres at 
closed points of the special fibre. We refer the reader to 
the exposition in [Pr] for a discussion of patching results 
and for detailled references on the subject. 
\par
Let $R$ be a complete discrete 
valuation ring with fractions field $K$, 
residue field $k$, and uniformiser $\pi$. Let $X$ be an 
admissible formal $R$-scheme which is an 
$R$-{\bf curve}, by which we mean that the special fibre 
$X_k:=X{\times}_R k$ is a 
reduced one dimensional $k$-scheme of finite type. Let $Z$ be a finite set of 
closed points of $X_k$. For a point $x\in Z$, let 
$X_x:=\Spf \hat O_{X,x}$ be the formal completion of $X$ at $x$, which is the
{\bf formal fibre} at the point $x$. Also let 
$X'$ be a formal open subset of $X$ whose special fibre is $X_k-Z$. For each 
point $x\in Z$, let $\{\wp_i\}_{i=1}^n$ be the set of
minimal prime ideals of $\hat {\Cal O}_{X,x}$ which contain $\pi$, they 
correspond to the {\it branches} $(\eta _i)_{i=1}^n$
of the completion of $X_k$ at $x$, and let $X_{x,i}:=\Spf 
\hat \Cal O_{x,\wp _i}$ be the formal completion of the localisation of $X_x$ 
at $\wp _i$. The ring $\hat \Cal O_{x,\wp _i}$ is a complete discrete 
valuation ring. The set 
$\{X_{x,i}\}_{i=1}^n$ is the set of {\bf boundaries} of the formal fibre
$X_x$. For each $i\in \{1,n\}$ we have a canonical morphism $X_{x,i}\to X_{x}$.

\pop {.5}
\par
\noindent
{\bf \gr 1.1. Definition.}\rm\ With the same notations as above a 
{\it (G)-cover patching data} for the pair $(X,Z)$ consists of the following:
\par
\noindent
a )\ A finite (Galois) cover $Y'\to X'$ (with group $G$).
\par
\noindent
b )\ For each point $x\in Z$, a finite (Galois) cover $Y_x\to X_x$ 
(with group $G$).

\par
The above datas a) and b) must satisfy to the following condition:

\par
\noindent
c )\ If $\{X_{x,i}\}_{i=1}^n$ are the boundaries of the formal fibre at the 
point $x$, then for each $i\in \{1,n\}$ is given a ($G$-equivariant) 
isomorphism $\sigma _i: Y_x\times _{X_x} X_{x,i}\simeq Y'\times _{X'} X_{x,i}$.

\pop {.5}
\par
\noindent
{\bf \gr 1.2. Proposition / Formal patching.}\rm\ {\sl Given a (G)-cover 
patching data as in 1.1 there exists a unique, up to isomorphism, 
(Galois) cover
$Y\to X$ (with group $G$) which induces the above (G)-covers 
in a) (resp. in b) when restricted to $X'$ (resp. when pulled back 
to $X_x$ for each point $x\in Z$).}

\pop {.5}
\par
The proof of 1.2 is an easy consequence of theorem 3.4. in [Pr], and which 
is due to Ferrand and Raynaud.

\pop {.5}
\par
\noindent
{\bf \gr 1.3. Remark.}\rm\ With the same notations as above let $x\in Z$
and let $\Tilde X_k$ be the normalisation of $X_k$. There is a one to one 
correspondance between the
set of points of $\Tilde X_k$ above $x$ and the set of boundaries of the 
formal fibre at the point $x$.
Let $x_i$ be the point of $\Tilde X_k$ above $x$
which corresponds to the boundary $X_{x,i}$, for $i\in \{1,n\}$. Assume 
that the point $x\in X_k(k)$ is rational. Then the completion of 
$\Tilde X_k$ at $x_i$ is isomorphic to the spectrum of a ring of formal 
power series
$k[[t_i]]$ in one variable over $k$, where $t_i$ is a local parameter at 
$x_i$. The complete local ring
$\hat \Cal O_{x,\wp _i}$ is a {\it two dimensional discrete valuation ring}
whose residue field is isomorphic to $k((t_i))$. Let $T_i$ be an element of 
$\hat \Cal O_{x,\wp _i}$ which lifts $t_i$, such an element is called a 
{\it parameter} of $\hat \Cal O_{x,\wp _i}$. Then it follows from [Bo] 
that there exists an isomorphism $\hat \Cal O_{x,\wp _i}\simeq R[[T_i]]
\{T_i^{-1}\}$, where $R[[T_i]]\{T_i^{-1}\}$ is the ring of formal power
 series $\sum_{i\in \Bbb Z} a_iT^{i}$ with $\lim _{i\to -\infty} \vert a_i
\vert =0$, and where $\vert\ \vert$ is an absolute value of $K$ associated 
to its valuation.
\pop {.5}
\par
\noindent
{\bf \gr 1.4. Rigid patching.} \rm  The analog of the above result 
is well known in rigid geometry, which is not surprising because of 
the link between rigid and formal geometry (cf. [Ra], [Bo-Lu]). We explain 
briefly the patching procedure in this context locally. 
Let $\Cal X:=\Spm A$ be an 
affinoid reduced curve, and let $X$ be a 
formal model of $\Cal X$ with special fibre $X_k$. Let $x$ be a 
closed point of $X_k$, and let $\Cal X_x$ be the {\it formal fibre} of 
$\Cal X$ at
$x$, which is a non quasi-compact rigid space and which consists of the 
set of points of $\Cal X$ which reduce to the point $x$. The structure 
of the {\it boundary} of $\Cal {X}_x$ is well known and depends 
functorially on the normalisation of the complete 
local ring $\hat {\Cal O}_{{X_k},x}$ (cf. [Bo-Lu-1]). Namely this boundary
decomposes 
into a disjoint union
of {\it semi-open annuli} $\Cal X_i:=\Cal X_{x,i}$ one corresponding to each 
minimal prime ideal $\eta _i$ of $\hat {\Cal O}_{X_k,x}$. Let 
$\Cal X':=\Cal X-\Cal 
X_x$ which is a quasi-compact rigid space. Let $f':\Cal Y'\to \Cal X'$ 
and $f_x:\Cal Y_x\to \Cal X_x$ be (Galois)-covers (with group $G$).
The (G)-cover $f':\Cal Y'\to \Cal X'$ extends to (G)-covers
$f_i:\Cal Y_i\to \Cal X_i$ of each of the components of the boundary of 
$\Cal X_x$, and the germ of such an extension is unique (cf. [Ra]). 
A (G)-patching 
data in this context are (G-equivariant) isomorphisms between the germs of 
$f_i:\Cal Y_i\to \Cal X_i$ and the restriction 
of the initial (G)-cover  
$f_x:\Cal Y_x\to \Cal X_x$ to $\Cal X_i$. The rigid patching result is that 
given a (G)-patching data as above then there exists a unique, up to 
isomorphism,
(G)-cover $f:\Cal Y\to \Cal X$ which induces the above covers above 
$\Cal X'$ and $\Cal X_x$ when restricted to this analytic subspaces.

\vskip.4cm
\epsfysize=5cm
\centerline{\epsfbox{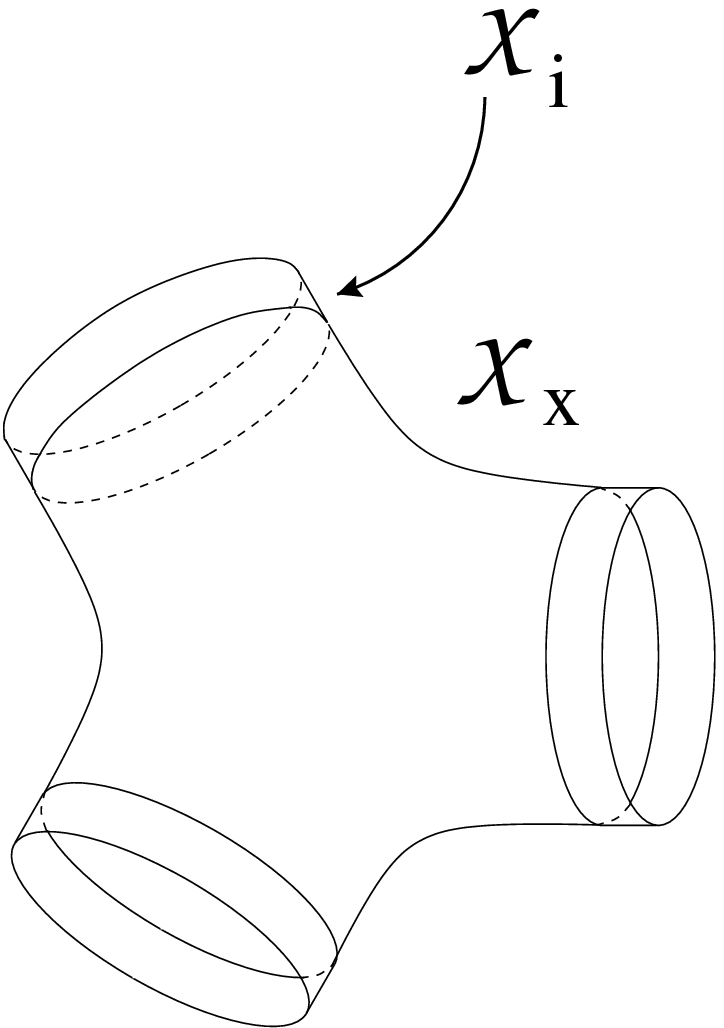}}
\vskip.4cm

\pop {1}
\par
\noindent
{\bf \gr  1.5. Local-global principle.}\rm\  As a direct consequence 
of the above patching results one obtains a {\it local-global principle}, 
which is certainly well known to the experts, for lifting of covers of 
curves. More precisely we have the following:

\pop {.5}
\par
\noindent
{\bf \gr  1.6. Proposition.}\rm\ {\sl let $X$ be a proper and flat algebraic 
(or formal) $R$-curve 
and let $Z:=\{x_i\}_{i=1}^n$
be a finite set of closed points of $X$. Let $f_k:Y_k\to X_k$ be a finite 
generically separable (Galois)-cover (of group $G$) whose branch locus is 
contained in $Z$. Assume that for each $i\in \{1,n\}$ there exists a 
(Galois)-cover $f_i:Y_{i}\to \Spf \hat {\Cal O}_{X,x_i}$ (of group $G$) 
which lifts the cover $\hat Y_{k,i}\to \Spec \hat {\Cal O}_{X_k,x_i}$ 
induced by $f_k$, where $\hat Y_{k,i}$ denotes the completion of $Y_k$ 
above $x_i$. Then 
there exists a unique, up to isomorphism, (Galois)-cover 
$f:Y\to X$ (of group $G$) which lifts the cover $f_k$, and which is 
isomorphic to the cover $f_i$ when pulled back to 
$\Spf \hat {\Cal O}_{X,x_i}$, for each $i\in \{1,n\}$}.
\pop {.5}
\par
\noindent
{\bf \gr Proof.}\rm\  After passing to the formal completion of $X$ along
its special fibre we may reduced to the case where $\Cal X$ is a formal 
$R$-curve. We treat the case where $Z=\{x\}$ consists of one point 
(the general case is similar). Let $U_k:=\Cal X_k-\{x\}$, and let 
$\Cal U$ be a formal 
open in $\Cal X$ whose special fibre equals
$U_k$. The \'etale cover $f'_k:V_k\to U_k$ induced by $f_k$ above $U_k$ 
can be lifted, uniquely, by the theorems of
lifting of \'etale covers (cf. [Gr]) to an \'etale formal cover $f':\Cal V\to 
\Cal U$. Let $\{\wp_i\}_{i=1}^n$ be the minimal prime ideals of 
$\hat {\Cal O}_{X,x}$ which contain $\pi$, and let $\Cal X_i:=\Cal X_{x,i}:
=\Spf \hat \Cal O_{x,\wp _i}$ be the formal completion of the localisation of 
$\hat {\Cal O}_{X,x}$ at $\wp _i$. We have canonical morphisms 
$\Cal X_i\to \Cal X$, and  $\Cal X_i\to
\Spf \hat \Cal O_{X,x}$. The cover $f'$ (resp the given cover
$f_i:Y_{i}\to \Spec \hat {\Cal O}_{X,x_i}$) induces (by pull back) 
a cover $f'_i:\Cal Y_i\to \Cal X_i$ (resp. a cover
$f_i:\Cal Y_i\to \Cal X_i$). For each $i\in \{1,n\}$, the cover $f_i$ and $f'_i$
by construction are isomorphic when restricted to the special fibre
$\Spec k((t_i))$ of $\Cal X_i$. Since both $f_i$ and $f'_i$ are \'etale 
and $\Cal X_i$ is local and complete we deduce that they are isomorphic.
Hence we obtain a patching data which allows us to patch the covers 
$f':\Cal V\to \Cal U$ and $f_i:Y_{i}\to \Spf \hat {\Cal O}_{X,x_i}$
in order to obtain a cover $f:\Cal Y\to \Cal X$ with the required properties.
Now thanks to the formal GAGA theorem this cover is algebraic $f:Y\to X$ 
and has the desired properties. Moreover if the starting datas are Galois 
then the constructed cover is also Galois with the same Galois group.

\pop {.5}
\par
\noindent
{\bf \gr  1.7. Remark.}\rm\ Although the formal patching result 1.2 and 
the rigid patching result 1.4 are equivalent, we opted in this paper to use
the formal patching result and the framework of formal geometry since
this seems to be more convenient for most readers. However it should be clear 
that one could also use the framework of rigid geometry and adapt the content
of this paper to this setting.

\pop {1}
\par
\noindent
{\bf \gr II. Computation of vanishing cycles.}
\pop {.5}
\par 
The main result of this section is 2.4 which gives a {\bf formula} 
which compares the dimensions of the space of vanishing cycles 
in a Galois cover 
$\Tilde f:\Cal Y\to \Cal X$ of group $\Bbb Z/p\Bbb Z$ between formal germs
of $R$-curves, where $R$ is a complete discrete valuation ring of inequal 
characteristic which contains a primitive $p$-th root of unity, where $p$ is 
the residue characteristic, in terms of the degeneration type of $\tilde f$ 
above the boundaries of $\Cal X$.

\pop {.5}
\par
\noindent
{\bf \gr 2.1.}\rm \ In this section we consider 
a complete discrete valuation ring $R$ of inequal 
characteristic, with residue characteristic $p>0$, and
which contains a primitive $p$-th root of unity $\zeta$. 
We denote by $K$ the fraction field of $R$, 
$\pi $ a uniformising parameter of $R$, $k$ the residue field of $R$, 
and $\lambda:=\zeta -1$. We also denote by $v_K$ the valuation of $K$ 
which is normalised by $v_K(\pi)=1$. {\bf We assume that the residue field
$k$ is algebraically closed}. By a (formal) $R$-curve we 
mean a (formal) $R$-scheme of finite type 
which is normal, flat, and whose fibres have dimension $1$. For an $R$-scheme 
$X$, we denote by $X_K:=X\times _{\Spec R}\Spec K$ the {\it generic} 
fibre of $X$, and $X_k:=X\times _{\Spec R}\Spec k$ its {\it special} fibre. 
In what follows by a {\it (formal) germ} $\Cal X$ 
of an $R$-curve we mean that
$\Cal X :=\Spec \Cal O_{X,x}$ is the (resp. $\Cal X:=\Spf \hat {\Cal O}_{X,x}$
is the formal completion of the) spectrum of 
the local ring of an
$R$-curve $X$ at a closed point $x$. Let $\Cal O_x$ be the local ring of 
$\Cal X_k$ at 
$x$. Let $\delta _x:=\dim _k \Tilde {\Cal O_x}/\Cal O_x$ where 
$\Tilde {\Cal O_x}$ is the normalisation of $\Cal O_x$ in 
its total ring of fractions, and let $r_x$ be the number of maximal ideals
in $\Tilde {\Cal O}_x$. The {\bf contribution to the arithmetic genus}
of the point $x$ is by definition $g_x:=\delta_x-r_x+1$. We will call the
integer $g_x$ the {\bf genus} of the point $x$. The following lemma 
is easy to prove (cf. for example [Bo-Lu-1]).

\pop {.5}
\par
\noindent
{\bf \gr 2.1.1. Lemma.}\rm\ {\sl Let $X_k$ be a proper reduced algebraic 
curve over $k$. Let $\tilde X_k\to X_k$ be the normalisation of $X_k$, and let
$\{X_i\}_{i\in I}$ be the irreducible components of $\tilde X_k$. Let
$\{x_j\}_{j\in J}$ be the singular point of $X_k$, which we assume to be 
rational. Let $g(X_k)$ (resp. $g(X_i)$) be the arithmetic genus of 
$X_k$ (resp. the arithmetic genus of $X_i$). Then 
$g(X_k)=\sum _{i\in I}g(X_i)+\sum _{j\in J}g_{x_{j}}$.}

\pop {.5}
\par
\noindent
{\bf \gr 2.2.}\rm\ Let $f:Y\to X$ be a finite cover between $R$-curves. 
Assume that the special fibres $X_k$ and $Y_k$ are reduced. Let $y$ 
be a closed point of
$Y$ and let $x$ be its image in $X$, which we assume to be a rational
point. Let $(x_j)_{j\in J}$ be the points of
the normalisation $\Tilde {X}_k$ of $X_k$ above $x$, and for a fixed
$j$ let $(y_{i,j})_{i\in I_j}$ be the points of the normalisation $\Tilde
{Y}_k$ of $Y_k$ which are above $x_j$. Assume that the morphism
$f_k: Y_k\to X_k$ is {\bf generically \'etale}. Under
this assumption we have the following {\it local Riemann-Hurwitz formula} 
which is due to Kato
(cf. [Ka], and [Ma-Yo]):

$$(1)\ (g_y+\delta_y-1)=n(g_x+\delta_x-1)+d_K-d^w_k$$
                                 
Where $n$ is the local degree at $y$, which is the degree of the morphism $
\Spec \hat \Cal O_{Y ,y} \to \Spec \hat \Cal O_{X,x}$ between the completion
of the local rings of $Y$ (resp $X$) at the point $y$ (resp. $x$), $d_K$
is the degree
of the divisor of ramification in the morphism
$\Spec (\hat \Cal O_{Y ,y}\otimes_R K) \to \Spec (\hat \Cal O_{X,x}\otimes_R K)
$.
Let $d^w_{i,j}:=v_{x_j}(\delta _{y_{i,j},x_j})-e_{i,j} +1$, where
$\delta _{y_{i,j},x_j}$
is the discriminant ideal of the extension $\hat \Cal O_{\Tilde {X}_k,x_j}
\to \hat \Cal O_{\Tilde {Y}_k,y_{i,j}}$ of complete discrete valuation
ring, and $e_{i,j}$   its ramification index. The integer $d^w_k$ is
equal to the sum
$\sum _{i,j}d^w_{i,j}$. In 2.4 we will obtain a formula similair to (1) in 
the case where $f$ is Galois of group $\Bbb Z/p\Bbb Z$ and which 
includes the case where $f_k$ is generically radicial.

\pop {.5}
\par
\noindent
{\bf \gr 2.3. Compactification process.}\rm \ Let $\Cal X 
:=\Spf \hat {\Cal O}_{X,x}$ be the formal germ of an $R$-curve at a 
closed point 
$x$. Let $\Tilde f:\Cal Y\to \Cal X$ be a Galois cover of group 
$\Bbb Z/p\Bbb Z$ with $\Cal Y$ local. We assume that the special fibre of 
$\Cal Y_k$ is reduced (this can alaways be achieved after a finite 
extension of $R$).
We will construct a compactification of the above cover $\Tilde f$ which 
will allows 
us to compute the arithmetic genus of the closed point of $\Cal Y$. 
More precisely we will construct a Galois cover $f:Y\to X$ of degree $p$,
between proper algebraic $R$-curves, a closed point $y\in Y$ and its 
image $x=f(y)$,
such that the formal germ of $X$ (resp. Y)
at $x$ (resp. at $y$) equals $\Cal X$ (resp. $\Cal Y$), and such that
the Galois cover 
$f_x:\Spf {\hat {\Cal O}}_{Y,y}\to \Spf {\hat {\Cal O}}_{X,x}$ induced by 
$f$ between the formal germs at 
$y$ and $x$ is isomorphic to the above given cover $\tilde f:\Cal Y\to 
\Cal X$. The
construction of such a compactification is well known in the case where
$\Tilde f$ has \'etale reduction type on the boundaries 
(cf [Ma-Yo] and [Ra-2]). In the case of radicial reduction type of degree 
$p$ on the boundaries one is able to 
carry out such a construction using the formal patching result in
1.2 and the result 3.1 in [Sa]. In fact it suffices 
to be able to treat the case of one boundary, which is easily done using 
the next proposition and the formal patching result.

\pop {.5}
\par
\noindent
{\bf \gr 2.3.1. Proposition.}\rm\ {\sl Let $D:=\Spf R<1/T>$ be the formal
closed disc centered at $\infty$. Let $\Cal D:=\Spf R[[T]]\{T^{-1}\}$, 
and let $\Cal D\to D$ be the canonical morphism. Let $\Tilde 
f:\Cal Y\to \Cal D$ be a non trivial
torsor under a finite and flat $R$-group scheme of rank $p$, such that
the special fibre of $\Cal Y$ is reduced. Then there exists
a Galois cover $f:Y\to D$ with group $\Bbb Z/p\Bbb Z$ whose pull back
to $\Cal D$ is isomorphic to the above given torsor $\Tilde f$. More precisely,
with the same notations as in [Sa] 3.1 we have the following depending on 
the several cases that occur there:
\par
case a ) consider the cover $f:Y\to D$ given generically by an equation 
$Z^p=\lambda ^pT^m+1$. This cover is an \'etale torsor above $D$ 
under the group scheme $\Cal H_{v_K(\lambda)}$ and induces an \'etale 
torsor $f_k:Y_k\to D_k$ in reduction. Moreover the genus of the 
smooth compactification  of $Y_k$ equals $(m-1)(p-1)/2$.
\par
case b-1 ) Consider the cover $f:Y\to D$ given generically by an equation
$Z^p=T^{-r}$ where $r:=p-h$. This cover is ramified at the generic fibre 
only above the point at infinity, the finite morphism $f_k:Y_k\to X_k$ is a 
$\mu_p$-torsor outside infinity, and $Y_k$ is smooth at infinity. Moreover 
the  genus of the smooth compactification of $Y_k$ equals $0$.
\par
case b-2 ) Consider the cover $f:Y\to D$ given generically by an equation
$Z^p=T^{-\alpha} (T^{-m}+1)$, where $\alpha$ is an integer such that $m+
\alpha \equiv 0 \mod\ (p)$. This cover is ramified at the generic fibre 
above infinity and above the distincts $m$-th root of unity. The finite 
morphism $f_k:Y_k\to D_k$ is a 
$\mu_p$-torsor outside infinity and the distincts $m$-th root of unity, 
and $Y_k$ is smooth. Moreover 
the  genus of the smooth compactification of $Y_k$ equals $0$.
\par
case c ) First if $m\le 0$ consider the cover $f:Y\to D$ given generically 
by an equation $Z^p=1+\pi ^{np}T^m$. This cover is a torsor under the
group scheme $\Cal H_n$, and its special fibre $f_k:Y_k\to X_k$ is a torsor 
under $\alpha _p$. The affine curve $Y_k$ has a unique singular point $y$
which is the point above infinity and $g_y=(-m-1)(p-1)/2$. Secondly if
$m\ge 0$ consider the cover $f:Y\to D$ given generically by an equation
$Z^p=T^{-\alpha}(T^{-m}+\pi ^{pn})$ where $\alpha$ is as in b-2. This
cover is ramified above infinity and the $m$ distinct roots of $\pi ^{-np}$.
The finite morphism $f_k:Y_k\to X_k$ is an $\alpha _p$-torsor outside 
infinity, and the special fibre
$Y_k$ of $Y$ is smooth. Moreover in both cases
the smooth compactification of the normalisation of $Y_k$ has genus $0$.}

\pop {.5}
\par
\noindent
{\bf \gr Proof.}\ Case a is straitforward. In case b-1 one has only to 
justify that $Y_k$ is smooth at infinity. the cover 
$f$ is given by the equation ${Z}^p=T^{-r}$, and after using the Bezout 
identity one reduces to an equation ${Z'}^p=T^{-1}$ from which it follows
directly that the complete local ring at the point above infinity is
$B=R[[Z']]$ hence this point is smooth. In case b-2 one deduces in a
similar way as above that 
$Y_k$ is smooth. Case c: if $m$ is negatif the $\alpha _p$-torsor
$f_k:Y_k\to X_k$ is given locally for the \'etale topology above the point
at infinity by an equation $z^p=t^{-m}$ where $t:=T\mod \ (\pi)$. The 
computation of the arithmetic 
genus of the singularity above infinity follows then by a direct calculation
(cf. example [Ra-1], or [Sai], 2.9). If $m$ is positive and in order to see that
$Y_k$ is smooth one considers the Galois cover $f':\Cal Y'\to \Cal P$ 
above the formal $R$-projective line $\Cal P$, obtained by gluing $D$ with 
the formal unit closed disc $D':=\Spf R<T>$ centered at $0$, and given by
the (same) equation $Z^p=1+\pi ^{np}T^m$. The genus of the generic fibre of 
$\Cal Y'$ is $(m-1)(p-1)/2$. The finite morphism $f'_k:\Cal Y'_k\to \Cal P_k$
is an $\alpha _p$-torsor outside infinity, and outside $0$ coincides with the
morphism $f_k$. Above $0$ this torsor is given locally for the \'etale topology
by an equation $z^p=t^m$ where $t:=T\mod \ (\pi)$, hence the arithmetic 
genus above $0$ equals $(m-1)(p-1)/2$ from which one deduces that $Y_k$
is smooth.

\pop {1}
\epsfysize=4cm
\centerline{\epsfbox{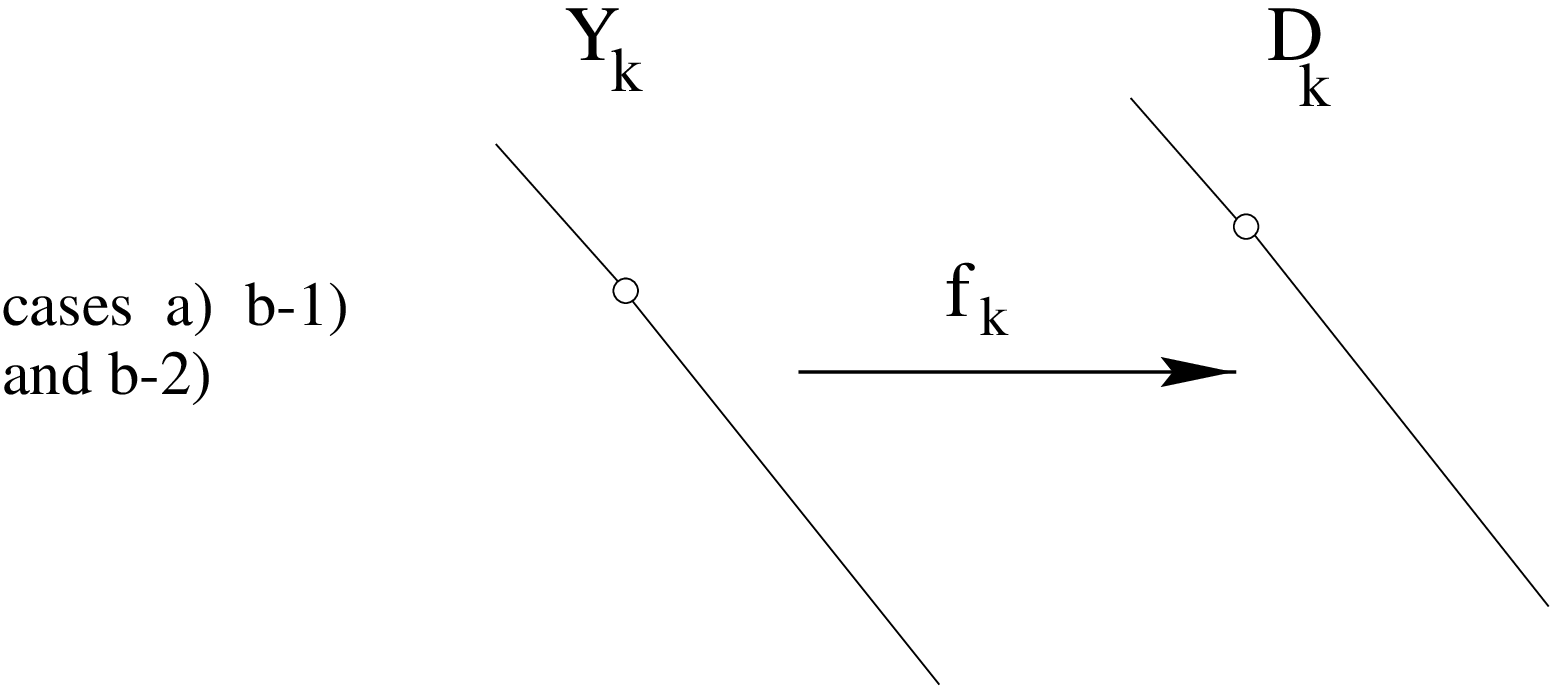}}

\pop {2}
\epsfysize=6cm
\centerline{\epsfbox{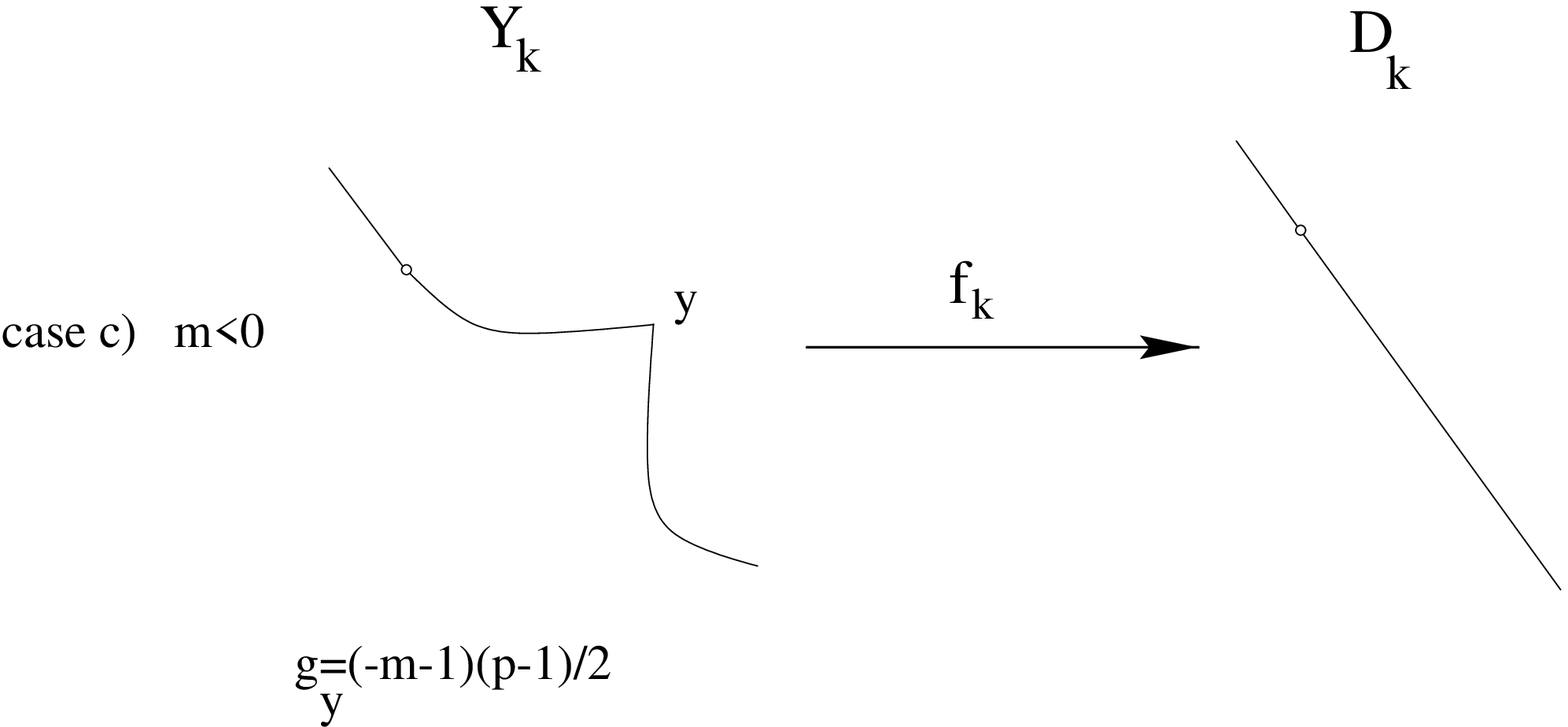}}

\pop {1}
\par
\noindent
{\bf \gr 2.3.2. Proposition}. \rm\ {\sl Let $\Cal X:=\Spf {\hat \Cal O_x}$ be 
the formal germ of an $R$-curve at a closed point $x$, and let $\{\Cal X_i\}
_i^n$
be the boundaries of $\Cal X$. Let 
$\tilde f:\Cal Y\to \Cal X$ be a 
Galois cover of group $\Bbb Z/p\Bbb Z$ with $\Cal Y$ local. Assume that 
$\Cal Y_k$ and $\Cal X_k$ are reduced. Then there 
exists a Galois cover $f:Y\to X$ of degree $p$,
between proper algebraic $R$-curves, a closed point $y\in Y$ and its 
image $x=f(y)$,
such that the formal germ of $X$ (resp. Y)
at $x$ (resp. at $y$) equals $\Cal X$ (resp. $\Cal Y$), and such that
the Galois cover 
$\Spf {\hat {\Cal O}}_{Y,y}\to \Spf {\hat {\Cal O}}_{X,x}$ induced by 
$f$ between the formal germs at 
$y$ and $x$ is isomorphic to the above given cover $\tilde f:\Cal Y\to 
\Cal X$. Moreover the formal completion of $X$ along its special fibre has a 
covering which consists of $n$ closed formal discs $
D_i$ which are patched
with $\Cal X$ along the boundaries $\Cal D_i$, and the special fibre $X_k$ 
of $X$ consists of $n$ smooth
projective lines which intersect at the point $x$.}

\pop {1}
\epsfysize=6.5cm
\centerline{\epsfbox{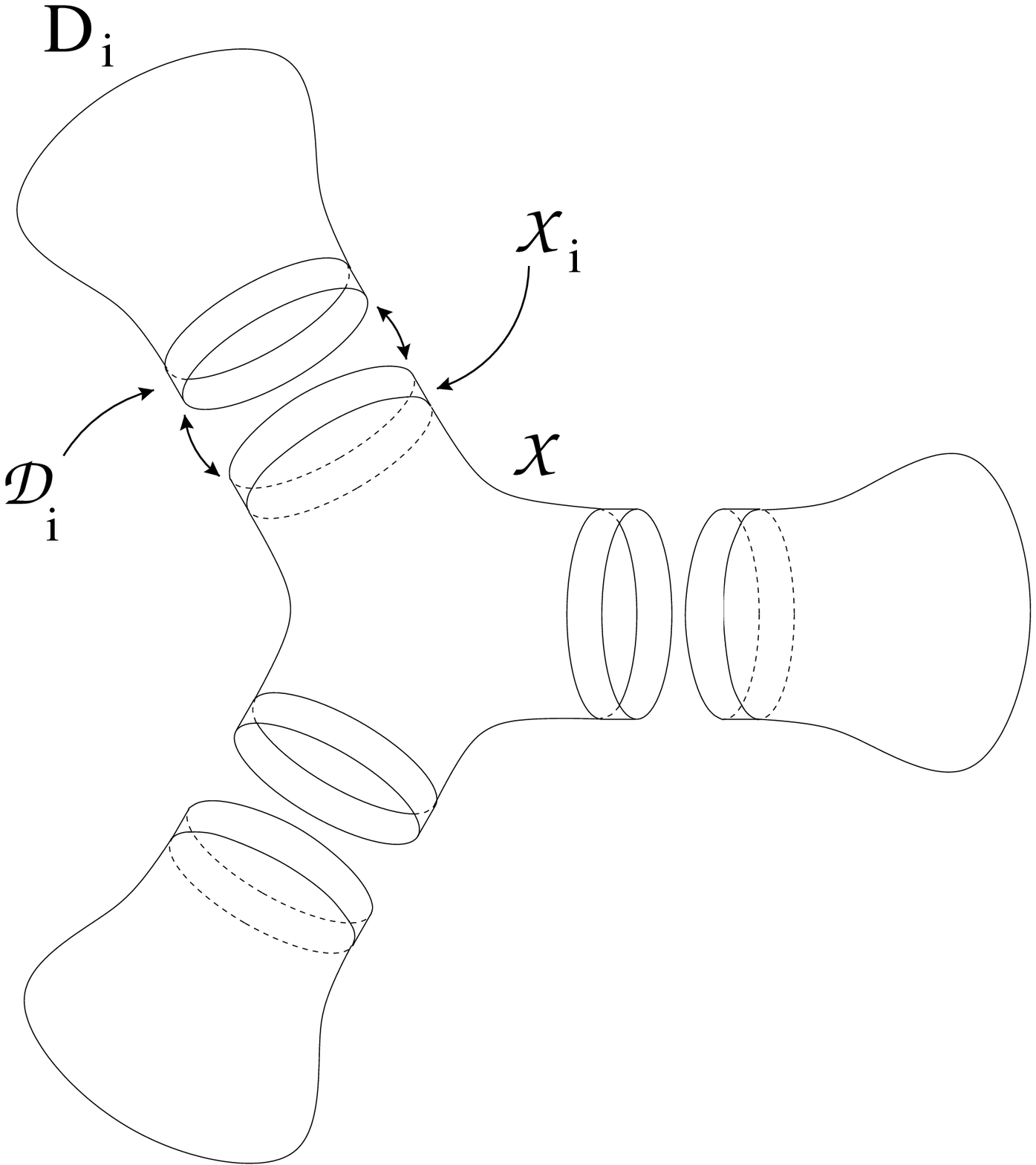}}

\pop {3}
\epsfysize=5cm
\centerline{\epsfbox{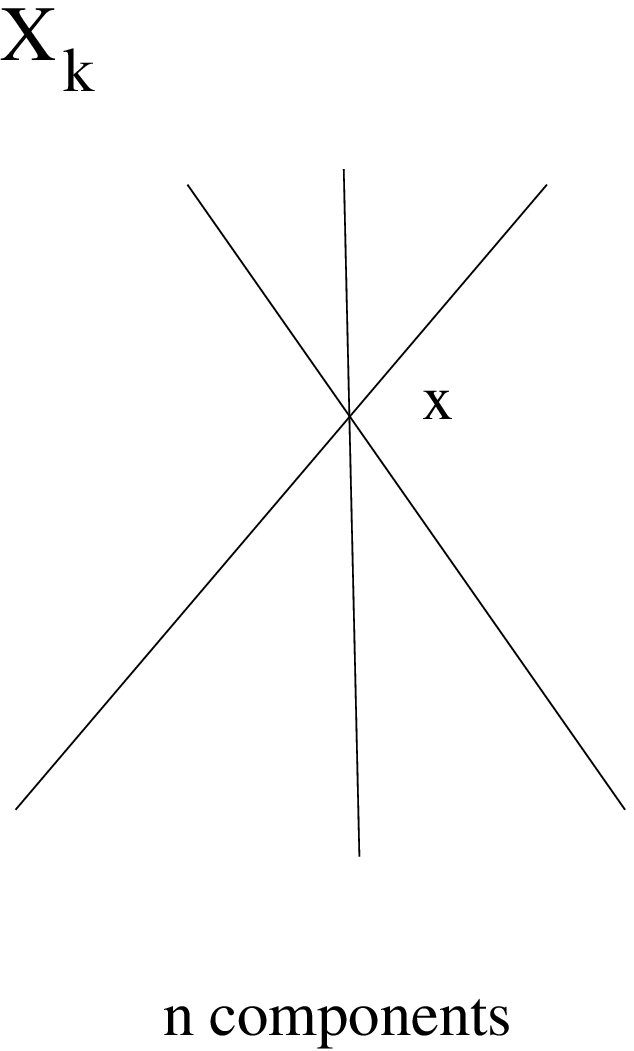}}

\pop {.5}
\par
\noindent
{\bf \gr Proof.}\rm \  Let $\{\wp_i\}_{i=1}^n$ be the minimal prime ideals of 
$\hat {\Cal O}_x$ which contain $\pi$ and which correspond to the  
branches $(\eta _i)_{i=1}^n$ of $\Cal X_k$ at $x$, and let $\Cal D_i:=\Spf 
\hat \Cal O_{\wp _i}$ be the formal completion of the localisation of $\Cal X$ 
at $\wp _i$. If $T_i$ is a 
lifting of a uniformising parameter of the branche $(\eta _i)$ 
of $\Cal X_k$ at $x$ then $\hat \Cal O_{\wp _i}$ is isomorphic to
$R[[T_i]]\{T_i^{-1}\}$. For 
each $i\in \{1,n\}$ consider a formal closed disc
$D_i:=\Spf R<1/T_i>$ centered at infinity, and the canonical morphism 
$\Cal D_i:=\Spf R[[T_i]]\{T_i^{-1}\}\to D_i$. As a consequence of the
formal patching result, which is valid for coherent sheaves 
(cf. [Pr] theorem 3.4), one can 
patch $\Cal X:=\Spf {\hat O}_x$
with the $D_i$, via the choice for each $i$ of an automorphism of $\Cal D_i$,
in order to construct a proper formal $R$-curve $X$, and a closed point 
$x\in X$, such that the formal completion of $X$ at $x$ equals $\Cal X$.
The special fibre $X_k$ of $X$ is a union of $n$ smooth $k$-projective 
lines which intersect at the point $x$.
Now the given cover $\tilde f$ induces a torsor 
$f_i:\Cal Y_i\to \Cal D_i$, under a finite and flat $R$-group scheme of 
rank $p$ for each $i$, and by the above lemma 3.3.1 one can 
find Galois covers $Y_i\to D_i$ of degree $p$ which after pull back to
$\Cal D_i$ coincide with $f_i$, for each $i\in \{1,n\}$. The formal patching
result again allows us then to patch these covers in order to construct a 
Galois cover $f:Y\to X$ of degree $p$ with the desired properties. The
formal $R$-curve is proper. By the formal GAGA theorems $X$ is algebraic 
and the Galois cover $f:Y\to X$ is also algebraic.

\vskip.4cm

\epsfysize=7cm
\centerline{\epsfbox{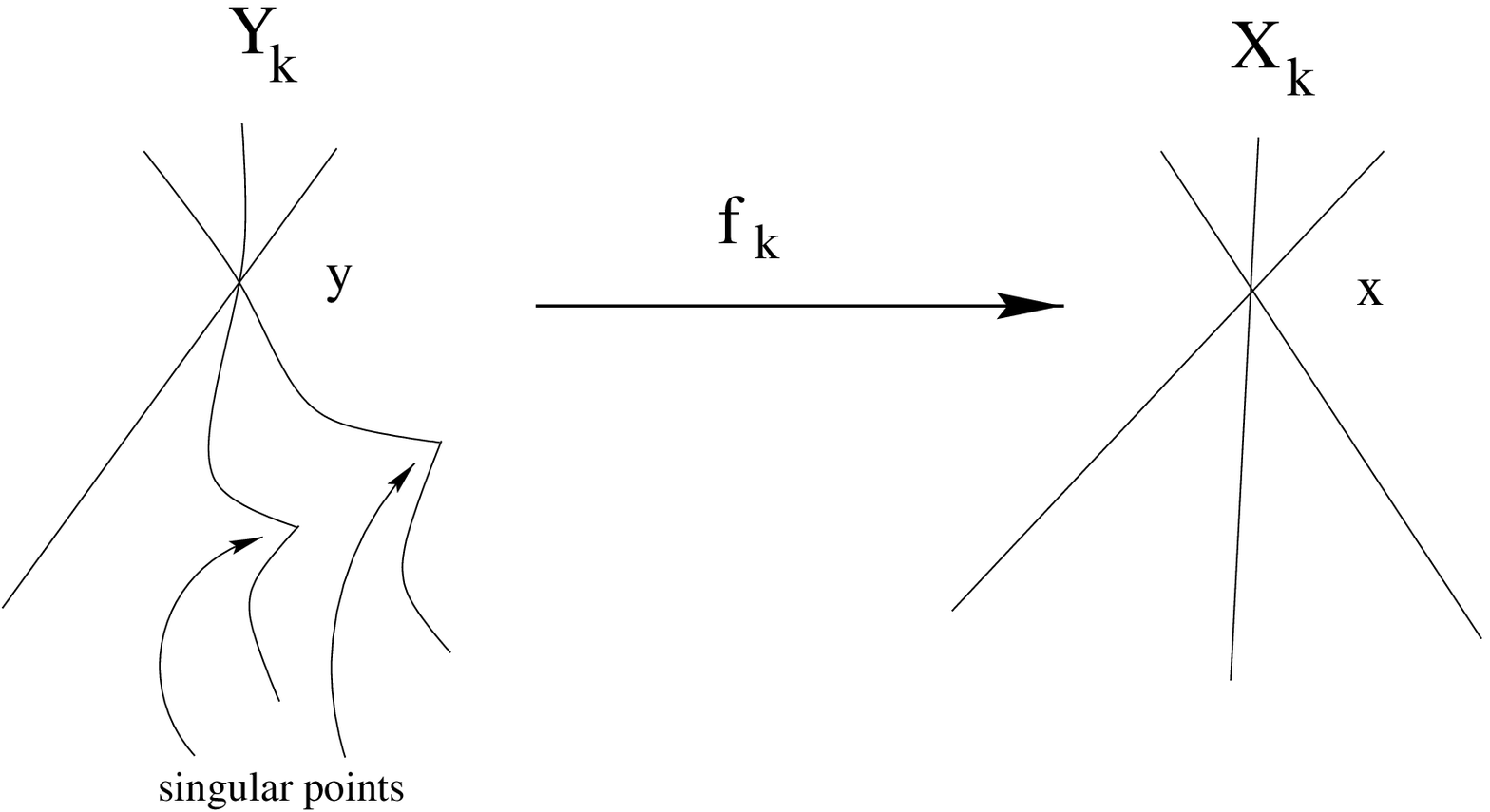}}

\pop {1}
\par
The next result is the main one of this paper, it gives a {\bf formula} which 
compares the dimensions of the space of vanishing cycles in a Galois 
cover of degree $p$ between formal fibres. 

\pop {.5}
\par
\noindent
{\bf \gr 2.4. Theorem.}\rm \ {\sl Let $\Cal X:=\Spf {\hat {\Cal O}_x}$ be 
the formal germ of an $R$-curve at a closed point $x$, with $\Cal X_k$ 
reduced. Let $\tilde f:\Cal Y\to \Cal X$ be a 
Galois cover of group $\Bbb Z/p\Bbb Z$ with $\Cal Y$ local, and 
$\Cal Y_k$ reduced. Let $\{\wp_i\}_{i\in I}$ be the minimal prime ideals of 
$\hat {{\Cal O}}_x$ which contain $\pi$ and which correspond to the  
branches $(\eta _i)_{i\in I}$ of $\Cal X_k$ at $x$, 
and let $\Cal X_{i}:=\Spf \hat \Cal O_{\wp _i}$ be the formal completion of 
the localisation of $\Cal X$ at $\wp _i$. For each $i\in I$, the above cover
$\tilde f$ induces a torsor $\tilde f_i:\Cal Y_i\to \Cal X_{i}$ under 
a finite and flat $R$-group scheme of rank $p$, above the boundary 
$\Cal X_i$. For 
each $i\in I$, let $(G_{k,i},m_i,h_i)$ be the
reduction type of $\tilde f_i$ (cf. [Sa] 3.2). Let $y$ be the closed point 
of $\Cal Y$. Then one has the following {\bf ``local Riemman-Hurwitz 
formula''}:
$$(2)\ \ \ 2g_y-2=p(2g_x-2)+d_{\eta}-d_s$$
Where $d_{\eta}$ is the degree of the divisor of ramification in the morphism
$\tilde f_K:\Cal Y_K\to \Cal X_K$ induced by $\tilde f$, where 
$\Cal X_K:=\Spec ({\hat {\Cal O}_x} \otimes _R K)$ and $\Cal Y_K:=
\Spec ({\hat {\Cal O}_{\Cal Y,y}}\otimes _R K)$, and $d_s:=
\sum _{i\in I^{\rad}}
(m_i-1)(p-1)+\sum _{i\in I^{\et}}(m_i-1)(p-1)$, where $I^{\rad}$ is the 
subset of $I$ consisting of those $i$ for which $G_{k,i}$ is radicial, and 
$I^{\et}$ is the subset of $I$ consisting of those $i$ for which 
$G_{k,i}$ is \'etale and $m_i\neq 0$.

\pop {.5}
\par
\noindent
{\bf \gr Proof.}\rm \ By 2.3.2 one can compactify the given morphism 
$\tilde f$. More precisely we constructed  a Galois cover $f:Y\to X$ of 
degree $p$ between proper algebraic $R$-curves, a closed point $y\in Y$ and 
its image $x=f(y)$, such that the formal germ of $X$ (resp. Y)
at $x$ (resp. at $y$) equals $\Cal X$ (resp. $\Cal Y$), and such that
the Galois cover $\Spf {\hat {\Cal O}}_{Y,y}\to \Spf {\hat {\Cal O}}_{X,x}$ 
induced by $f$ between the formal germs at $y$ and $x$ is isomorphic 
to the given cover $\tilde f:\Cal Y\to \Cal X$. The special fibre of $\Cal X$
consists (by construction) of $n$-ditincts smooth projective lines which
meet at the closed point $x$. The formal completion of $X$ along its 
special fibre has a covering
which consists of $\vert I\vert$ formal closed unit discs which are patched 
with the formal fibre $\Cal X$ along the boundaries $\Cal X_i$. 
The above formula 
(2) follows then by comparing the arithmetic genus of the generic fibre 
$Y_K$ of $Y$ and the one of the special fibre $Y_k$. By the precise 
informations given in
3.3.1 one can easily deduce that $g(Y_K)=p g_x+(1-p)+d_{\eta}/2+\sum 
_{i\in I_b} (-m_i+1)(p-1)/2+\sum _{i\in I_{c,>}} (-m_i+1)(p-1)$ where
$I_b$ is the subset of $I$ consisting of those $i$ for which the degeneration
data correspond to the one in case b of [Sa] 3.1, and $I_{c,>}$ is the subset 
of $I$ consisting of those $i$ for which the degeneration
data correspond to the one in case c in [Sa] 3.1, and with $m$ positif. On the 
other hand one has $g(Y_k)=g_y+\sum _{i\in I_{a,>}}(-m_i-1)(p-1)/2 +
\sum _{i\in I_{c,<}} (-m_i+1)(p-1)$, where 
$I_{a,>}$ is the subset 
of $I$ consisting of those $i$ for which the degeneration
data correspond to the one in case a of 2.4.1 with $m_i\neq 0$, and
$I_{c,<}$ is the subset 
of $I$ consisting of those $i$ for which the degeneration
data correspond to the one in case c in [Sa] 3.1, and with $m_i<0$. 
Now since $Y$ is flat $g(Y_K)=g(Y_k)$
from which directly follows the formula (2).

\pop {.5}
\par
\noindent
{\bf \gr III. Galois covers of degree $p$ above germs of semi-stable 
curves.}\rm \ 
In what follows we use the same notations as in 2.1. 
As a consequence of the above result 2.4 we will 
deduce some interesting results in this section 
in the case of a Galois cover $\Cal Y\to \Cal X$, where $\Cal X$ is the 
formal germ of a semi-stable $R$-curve at a closed point. These results 
will play an important role in 
the paper [Sa-1] in order to exhibit, and realise, the degeneration datas 
which describe the semi-stable reduction of Galois covres of degree $p$. 

\pop {.5}
\par
\noindent
{\bf \gr 3.1.}\rm\ We start with the case of 
a Galois cover of degree $p$ above a germ of a {\bf smooth} point.

\pop {.5}
\par
\noindent
{\bf \gr 3.1.1. Proposition.}\rm\ {\sl Let $\Cal X:=\Spf R[[T]]$
be the germ of a formal $R$-curve at a smooth point $x$. Let $\eta $ be the 
generic point of the special fibre $\Cal X_k$ of $\Cal X$. The completion 
of the localisation of $\Cal X$ at $\eta$ is $\Cal X_{\eta}:=
\Spf R[[T]]\{T^{-1}\}$, which is the boundary of $\Cal X$. Let 
$f:\Cal Y\to \Cal X$ be a Galois cover of degree $p$, with $\Cal Y$ 
local. 
Assume that the special fibre of $\Cal Y$ is 
reduced. Let $y$ be the unique closed point 
of $\Cal Y_k$. 
Let $\delta _K:=r(p-1)$ be the degree of the divisor of 
ramification in the morphism $f:\Cal Y_K\to \Cal X_K$. We 
distinguish two cases:
\par
\noindent
case 1 ) $\Cal Y_k$ is {\bf unibranche} at $y$. Let $(G_k,m,h)$ be the 
degeneration type of $f$ above the boundary $\Cal X_{\eta}$. Then
necessarily $r-m-1\ge 0$, and $g_y=(r-m-1)(p-1)/2$.
\par
\noindent
case 2 ) $\Cal Y_k$ has $p$-{\bf branches} at $y$. Then the cover $f$ has 
an \'etale completely split reduction of type $(\Bbb Z/p\Bbb Z,0,0)$  
on the boundary, i.e. the induced 
torsor above $\Spf R[[T]]\{T^{-1}\}$ is trivial, in which
case $g_y=(r-2)(p-1)/2$.}

\pop {2}

\epsfysize=8cm
\centerline{\epsfbox{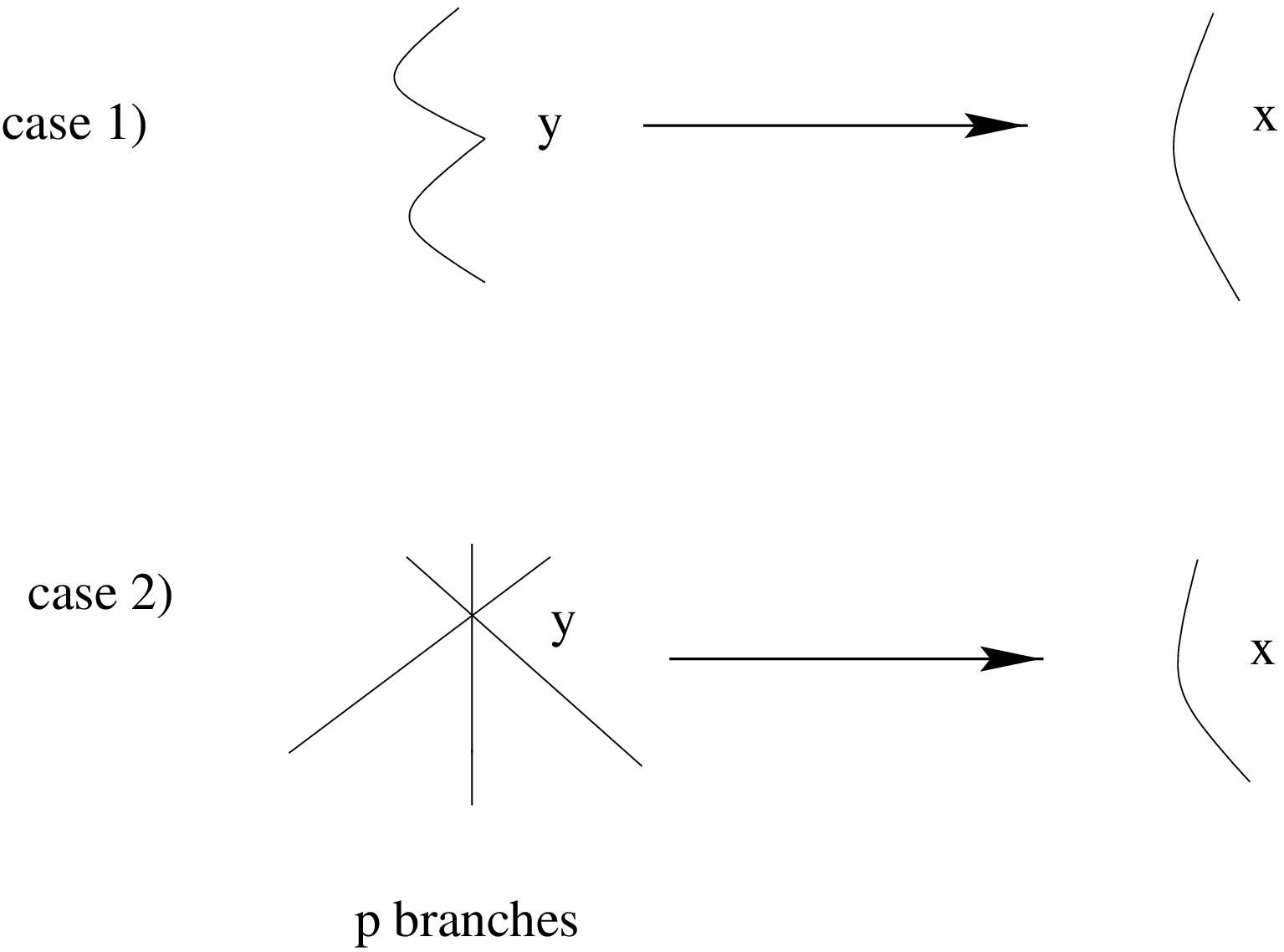}}

\pop {1}
\par
With the same notations as in 3.1.1, and as an immediate consequence, one can
immediately see whether the point $y$ is smooth or not. More precisely we 
have the following:

\pop {.5}
\par
\noindent
{\bf \gr 3.1.2. Corollary.}\rm\ {\sl We use the same notations as in 3.1.1 
Then $y$
is a smooth point, which is equivalent to $g_y=0$, if and
only if $r=m+1$. This in the case of radicial reduction type on the 
boundary is equivalent to $r=-\ord (\omega)$, where $\omega$ is the associated 
differential form. In particular if the reduction is of multiplicative type 
on the boundary, i.e. $G_k=\mu_p$, then $g_y=0$ only if $r=1$ or $r=0$, since 
$\ord (\omega)\ge -1$ in this case. Also if $r=1$ and $g_y=0$ then 
necessarily $G_k=\mu _p$}.

\pop {.5}
\par
Next we give some examples of Galois covers of degree $p$
above the formal germ of a smooth point which cover all the possibilities
for the genus and the degeneration type on the boundary. Both in 3.1.3 and 
3.1.4 we use the same notations as in 3.1.1. We first begin with        
examples with genus $0$.

\pop {.5}
\par
\noindent
{\bf \gr 3.1.3. Examples.}\rm\ The following are examples 
given by explicit equations of the different cases, depending on the 
possible degeneration type, of 
Galois covers $f:\Cal Y\to \Cal X$ of degree $p$ 
above $\Cal X=\Spf R[[T]]$, and where $\bold {g_y=0}$. 
\par
1 ) For $m>0$ an integer prime to $p$, consider the cover given generically 
by the equation $X^p=1+\lambda ^pT^{-m}$. 
Hier $r=m+1$, and this cover has a reduction of type 
$(\Bbb Z/p\Bbb Z,m,0)$ on the boundary.
\par
2 ) For $h\in \Bbb F_p^*$, consider the cover given generically 
by the equation $X^p=T^h$. Hier $r=1$, and this cover has a reduction of 
type $(\mu _p,0,h)$ on the boundary. 

\par
3 ) Consider the cover given generically 
by the equation $X^p=1+T$. Hier $r=0$, and this cover has a reduction of type
$(\mu _p,-1,0)$ on the boudary.

\par
4 ) For $n< v_K(\lambda)$, and $m<0$, consider the cover given generically 
by the equation $X^p=1+\pi  ^{np}T^m$. Hier $r=-m+1$, and this cover has a 
reduction of type $(\alpha _p,-m,0)$ on the boudary. 

\par
5) For $n< v_K(\lambda)$, consider the cover given generically 
by the equation $X^p=1+\pi  ^{np}T$. Hier $r=0$, and this cover has a 
reduction of type
$(\alpha _p,-1,0)$ on the boudary. 
\pop {.5}
\par
Next we give examples of Galois covers of degree $p$ 
above formal germs of smooth points which lead to a singularity with 
positive genus.

\pop {.5}
\par
\noindent
{\bf \gr 3.1.4. Examples.}\rm\ The following are examples 
given by explicit equations of the different cases, depending on the 
possible reduction type, of 
Galois covers $f:\Cal Y\to \Cal X$ of degree $p$ 
above $\Cal X=\Spf R[[T]]$, and where $\bold {g_y>0}$. 

\par
1 ) For $m>0$ an integer prime to $p$, and $m'>m$, consider the cover 
given generically 
by the equation $X^p=1+\lambda ^p(T^{-m}+\pi T^{-m'})$. 
Hier $r=m'+1$, and this cover has a reduction of type 
$(\Bbb Z/p\Bbb Z,m,0)$ on the boundary. Moreover the genus $g_y$ of 
the closed point $y$ of $\Cal Y$ equals $(m'-m)(p-1)/2$.
\par
2 ) For $h\in \Bbb F_p^*$, and $m>0$ an integer prime to $p$, 
consider the cover given generically 
by the equation $X^p=T^{h'}(T^m+a)$, where $h'$ is a positif integer such that 
$m+h'\equiv h\mod \ (p)$, and $a\in \pi R$. Hier $r=m+1$, and this cover 
has a reduction of type $(\mu _p,0,h)$ on the boundary. Moreover the genus
$g_y$ of the closed point $y$ of $\Cal Y$ equals $m(p-1)/2$.

\par
3 ) For a positif integer $m'$, an integer $h$ such that  
$m'+h\equiv 0\mod \ (p)$, and $a\in \pi R$, consider the cover given 
generically by the equation $X^p=T^h(T^{m'}+a)(1+T^m)$. Hier $r=m'+1$, and 
this cover has a reduction of type $(\mu _p,-m,0)$ on the boudary. Moreover 
the genus $g_y$ of the closed point $y$ of $\Cal Y$ equals $(m'+m)(p-1)/2$.

\par
4 ) For $n< v_K(\lambda)$, and integers $m>0$ prime to $p$ and $m<m'$, 
consider the cover given generically 
by the equation $X^p=1+\pi  ^{np}(T^{-m}+\pi T^{-m'})$. Hier $r=m'+1$, and 
this cover has a reduction of type $(\alpha _p,-m,0)$ on the boudary. Moreover 
the genus $g_y$ of the closed point $y$ of $\Cal Y$ equals $(m'+m)(p-1)/2$.

\pop {.5}
\par
Note that in cases 5, 6 of 3.1.3, and 4 of 3.1.4, and in order to realise 
these covers above $\Cal X$, for a given $n$, one needs in general to 
perform a ramified extension of $R$.

\pop {.5}
\par
\noindent
{\bf \gr 3.2.}\ \rm
Next we examine the case of Galois covers of degree $p$ above formal germs
at {\bf double points}.

\pop {.5}
\par
\noindent
{\bf \gr 3.2.1. Proposition.}\rm\ {\sl Let $\Cal X:=\Spf R[[S,T]]/(ST-\pi ^e)$
be the formal germ of an $R$-curve at an ordinary double point $x$ of 
thikeness $e$. The special fibre of 
$\Cal X$ consists of two irreducible components $X_1$ and $X_2$ with 
generic points $\eta _1$ and $\eta _2$, corresponding to the prime ideal $(\pi,
T)$ (resp. $(\pi ,S)$) of $R[[S,T]]/(ST-\pi ^e)$. The completion 
of the localisation of $\Cal X$ at $\eta_1$ (resp. $\eta _2$) is isomorphic 
to $\Cal X_1:=\Spf R[[S]]\{S^{-1}\}$ (resp. $\Cal X_2:=\Spf 
R[[T]]\{T^{-1}\}$). These
are the two boundaries of $\Cal X$. Let 
$f:\Cal Y\to \Cal X$ be a Galois cover of group $\Bbb Z/p\Bbb Z$ with 
$\Cal Y$ local. Assume that the special fibre of $\Cal Y$ is 
reduced. We assume that $\Cal Y_k$ has {\bf two branches} at the point $y$. 
Let $\delta _K:=r(p-1)$ be the degree of the divisor of 
ramification in the morphism $f:Y_K\to X_K$. Let $(G_{k,i},m_i,h_i)$
be the type of reduction on the two boundaries of $\Cal X$, for $i=1,2$. Then
necessarily $r-m_1-m_2\ge 0$, and $g_y=(r-m_1-m_2)(p-1)/2$}.

\pop {1}
\epsfysize=2.6cm
\centerline{\epsfbox{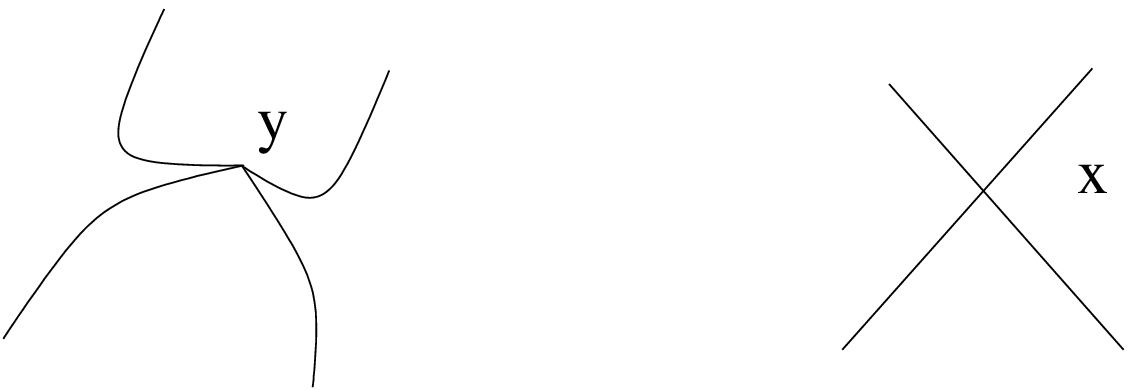}}

\pop {1}
\par
\noindent
{\bf \gr 3.2.2. Proposition.}\rm\ {\sl We use the same notations as in 
3.2.1. We consider the remaining cases:
\par
\noindent
case 1)
{\bf $\Cal Y_k$ has $p+1$ branches at} $y$. 
We can assume that $\Cal Y$ is completely split above $\eta _1$. 
Let $(G_{k,2},m_2,h_2)$ be the reduction type on the second boundary of 
$\Cal X$. Then necessarily $r-m_2-1\ge 0$ and $g_y=(r-m_2-1)(p-1)/2$.
\par
\noindent
case 2) {\bf $\Cal Y_k$ has $2p$ branches at $y$}. In this case $\Cal Y$ is 
completely split above $\eta _1$ and $\eta _2$ and $g_y=(r-2)(p-2)/2$.}

\pop {.5}
\par
With the same notations as in 3.2.1, and as an immediate consequence, one can
recognise whether the point $y$ is a double point or not. 
More precisely we have the following.

\pop {.5}
\par
\noindent
{\bf \gr 3.2.3. Corollary.}\rm\ {\sl We use the same notations as in 3.2.1. 
Then $y$ is an ordinary double point, which is equivalent to $g_y=0$, if and
only if $x$ is an ordinary double point of thikeness divisible by $p$, and
$r=m_1+m_2$. Moreover if $g_y=0$, if $r=0$, and if 
$(G_{k,i},m_i,h_i)$ is the reduction type on the boundary for $i=1,2$, then
necessarily $h_1+h_2=0$}.

\pop {.5}
\par
\noindent
{\bf \gr Proof.}\rm \ We need only to justify the last assertion. If both 
$h_1$ and $h_2$ equals $0$ there is nothing to prove. Otherwise assume 
$h_1\neq0$. Then $m_1=m_2=0$ (because $g_y=0=m_1+m_2$), and one sees easily 
that necessarily
$G_{k,1}=G_{k,2}=\mu_p$. So the cover $f:\Cal Y\to \Cal X$ is 
in this case a $\mu_p$-torsor and its 
reduction $f_k:\Cal Y_k\to \Cal X_k$ is also $\mu_p$-torsor given by an 
equation $t^p=u$, and $\omega :=du/u$ is the associated differential form. The 
restriction $\omega _i$ of $\omega$ to the $i$-th branch , for $i=1,2$, is 
the differential 
form associated to the $\mu_p$-torsor $f_i:\Cal Y_i\to \Cal X_i$ induced by 
$f$ above the boundary $\Cal X_i$ of $\Cal X$.
The equality  $h_1+h_2=0$ follows then from the fact that the sum of the 
residues of $\omega$ on each branch equals $0$, which is a property of the 
regular differential form at a double point.

\pop {.5}
\par
Next we give some examples of Galois covers of degree $p$,
above the formal germ of a double point, which leads to singularities with genus $0$, i.e. double points, and such that $r=0$. These examples will be used in
the paper [Sa-1] in order to realise the ``degeneration datas''.

\pop {.5}
\par
\noindent
{\bf \gr 3.2.4. Examples.}\rm \ The following are examples, given by 
explicit equations, of the different cases, depending on the possible type of 
reduction on the boundaries, of Galois covers $f:\Cal Y\to \Cal X$ of 
degree $p$ above $\Cal X=\Spf R[[S,T]]/(ST-\pi ^e)$, with $\bold {r=0}$, 
and where $\bold {g_y=0}$, for a suitable choice of $e$ and $R$. Note that 
$e=pt$ must be divisible by $p$. In all the following examples we have $r=0$.

\par
1 ) $p$-{\bf Purity}: if $f$ as above has an \'etale reduction type on the 
boundaries, and $r=0$, then $f$ is necessarily \'etale, and hence is 
completely split since $\Cal X$ is strictly henselian.
\par
2 ) Consider the cover given generically by an equation $X^p=T^h$, which 
leads to a reduction on the 
boundaries of type $(\mu_p,0,h)$ and $(\mu_p,0,-h)$.

\par
3 ) For a fixed integer $m>0$ prime to $p$, and after eventually a ramified 
extension of $R$ choose $t$ such that $tm< v_K(\lambda)$, and
consider the cover given 
generically by an equation $X^p=1+S^{m}$, which leads to a reduction on 
the boundaries of type $(\mu_p,-m,0)$ and $(\alpha_p,m,0)$.

\par
4 )  For a fixed integer $m>0$ prime to $p$, and after eventually a ramified 
extension of 
$R$ choose $t$ such that $t=v(\lambda)/m$, and consider the cover given 
generically by an equation 
$X^p=\lambda ^p/T^m+1$, which leads to a reduction on the boundaries
of type $(\Bbb Z/p\Bbb Z,m,0)$ and $(\mu_p,-m,0)$.

\par
5 ) For a fixed integer $m>0$ prime to $p$, and after eventually a 
ramified extension of $R$, choose $t$ such that $tm< v_K(\lambda)$, and 
consider the cover given generically by an equation
$X^p=1+\lambda ^pS^{-m}$, which leads to a reduction on the boundaries
of type $(\Bbb Z/p\Bbb Z,m,0)$ and $(\alpha_p,-m,0)$.
\par
6 ) For a fixed integer $m>0$ prime to $p$, and after eventually a ramified 
extension of $R$ choose $t$ and $n$ such that $tm+n< v_K(\lambda)$, 
and consider the cover given 
generically by an equation $X^p=1+\pi ^{np}S^m$, which leads to a reduction 
on the boundaries of type $(\alpha_p,-m,0)$ and $(\alpha_p,m,0)$.

\pop {.5}
\par
In fact one can describe Galois covers of degree $p$ above formal germs
at double points, which are \'etale above the generic fibre, and with
genus $0$. Namely they are all of the form given in the examples 3.2.4. 
In particular these covers are uniquely determined, up to isomorphism, by 
their degeneration type on the boundaries. More precisely we have the 
following:

\pop {.5}
\par
\noindent
{\bf \gr 3.2.5. Proposition.}\rm\ {\sl Let $\Cal X$ be the formal germ of a
semi-stable $R$-curve at an ordinary double point $x$. Let $f:\Cal Y\to \Cal X$
be a Galois cover of degree $p$, with $\Cal {Y}_k$ reduced and local, and 
with $f_K:\Cal Y_K\to \Cal X_K$ {\bf \'etale}. Let ${\Cal X_i}$, for $i=1,2$, 
be the boundaries of $\Cal X$. Let $f_i:\Cal Y_i\to \Cal X_i$ be the 
torsors induced by $f$ above $\Cal X_i$, and let $\delta _i$ be the 
corresponding degree of the different (cf. 2.3.5). 
Let $y$ be the closed point of $\Cal Y$. Assume that $\bold {g_y=0}$. 
Then there exists an isomorphism $\Cal X\simeq \Spf R[[S,T]]/(ST-\pi ^{tp})$ 
such that, if say $\Cal X_2$ is the boundary corresponding to the prime 
ideal $(\pi,S)$, the following holds:

\par
a )\ The cover $f$ is generically given by an equation $X^p=T^h$, with 
$h\in \Bbb F_p^*$,
which leads to a reduction on the boundaries of $\Cal X$ of type 
$(\mu_p,m=0,h)$ and 
$(\mu_p,m=
0,-h)$. Hier $t>0$ can be any integer. In this case $\delta _1
=\delta _2=v_K(p)$.
\par
b )\ The cover $f$ is generically given by an equation $X^p=1+T^m$, 
for $m>0$ an integer prime to $p$, such that $tm<v_K(\lambda)$. In 
particular in this case
necessarily $t<v_K(\lambda)$. This cover leads 
to a reduction on the boundaries of $\Cal X$ of type $(\alpha_p,m,0)$
and $(\mu_p,-m,0)$. In this case $\delta _2=v_K(p)=\delta _1+(p-1)tm$
\par 
c )\ The cover $f$ is generically given by an equation 
$X^p=1+T^m$, with $m>0$ an integer prime to $p$, such that 
$tm=v_K(\lambda)$. In particular in this case $t$ divides $v_K(\lambda)$.
This cover leads to a reduction on the boundaries of $\Cal X$
of type $(\Bbb Z/p\Bbb Z,m,0)$ and $(\mu_p,-m,0)$. In this case 
$\delta _2=v_K(p)=\delta _1+(p-1)tm$, and $\delta _1=0$.
\par
d )\ The cover $f$ is generically given by an equation
$X^p=1+\lambda ^pT^{-m}$, with $m>0$ an integer prime to $p$ 
such that $tm<v_K(\lambda)$, which leads to a reduction on the boundaries
of $\Cal X$ of type $(\Bbb Z/p\Bbb Z,m,0)$ and $(\alpha_p,-m,0)$. In 
particular we have necessarily $t<v_K(\lambda)$. In this case 
$\delta _1=\delta _2+(p-1)tm$, and $\delta _2=0$.
\par
e )\ The cover $f$ is generically given by an equation 
$X^p=1+\pi ^{np}T^m$, for a positif integer $m$ prime to $p$, with 
$n<v_K(\lambda)$, and $n+tm<v_K(\lambda)$, 
which leads to a reduction on the boundaries of type $(\alpha_p,-m,0)$ and 
$(\alpha_p,m,0)$}. In particular we have necessarily $t<v_K(\lambda)$.
In this case $\delta _2=\delta _1+(p-1)tm$, and 
$\delta _1=(p-1)(v_K(\lambda)-(n+tm))$.
\par
In all cases we have $\delta _2-\delta _1=mt(p-1)$.}

\pop {.5}
\par
\noindent
{\bf \gr Proof.}\rm \ We explain briefly the proof. Since $g_y=0$, the 
thikness $e=pt$ at the double point $x$ of $\Cal X$ is divisible 
by $p$. Say $\Cal X\simeq \Spf A'$, where 
$A'=R[[T',S']]/(S'T'-\pi ^{tp})$. The $\mu _p$-torsor 
$f_K:\Cal Y_K\to \Cal X_K$ is given by an equation $X^p=u_K$ where,
$u_K$ is a unit on $\Cal X_K$. Such a unit can be uniquely written as
$u_K:=\pi ^n {T'}^{m}u$, where $n$ and $m$ are integers, and $u\in
R[[T',S']]/(S'T'-\pi ^{tp})$ is a unit. Note first that necessarily 
$n\equiv 0\mod (p)$,
since $\Cal Y_k$ is reduced. Assume first that $\gcd (m,p)=1$. Let
$t':=T'\mod (\pi)$, let $s':=S'\mod (\pi)$, and let
$\bar u:=u\mod (\pi)$ which is a unit of $k[[s',t']]/(s't')$. 
Let $T:=T'u^{1/m}$. Then $\Cal X\simeq \Spf A$, where 
$A=R[[S,T]]/(ST-\pi ^{tp})$ for a suitable $S\in A$.
Let $B:=A[X,Y]/(X^p-T^m,Y^p-S^{m},XY-\pi^{t})$. Then $B$ is a finite flat 
algebra over $A$ which is integrally closed 
(because $B/\pi B$ is reduced), and $\Cal Y=\Spf B$. The cover 
$f:\Cal Y\to \Cal X$ is thus generically given by an equation $X^p=T^m$, and 
we are in case a. Assume now that 
$m\equiv 0 \mod (p)$. Hier we have two cases:
\par
case 1)\ $\bar u$ is not a $p$-power. Then it is easy to see, that after
changing the Kummer generator of the torsor $f_K$, one can assume that
$u$ is such that $\bar u=1+{t'}^m\bar v$, where $m$ is a positif integer 
prime to $p$ and $\bar v$ is a unit of $k[[s',t']]/(s't')$. In particular 
$\bar u=1+(t'\bar v ^{1/m})^m=1+t^m$, where 
$t:=t'\bar v^{1/m}$. Let $T:=T'v^{1/m}$, where $v:=(u-1)/{T'}^m$. 
We have $\Cal X\simeq \Spf R[[S,T]]/(ST-\pi ^{tp})$ for 
a suitable $S\in A$.  The cover $f:\Cal Y\to \Cal X$ is thus generically 
given by an equation $X^p=1+T^m$, and 
we are in case b. Let $\Cal X_1:=\Spf R[[S]]\{S^{-1}\}$ be the boundary of 
$\Cal X$ corresponding to the ideal $(\pi, T)$. The torsor
$f_1:\Cal Y_1\to \Cal X_1$ induced by $f$ above $X_1$ is generically 
given by an equation $X^p=1+\pi ^{pt}S^{-m}$, which imply that
$t\le v_K(\lambda)$, since $f_1$ is not completely split. Moreover
we are in case $b$ if $t<v_K(\lambda)$, and in case $c$ if $t=v_K(\lambda)$.
\par
case 2)\ $\bar u$ is a $p$-power. In this case, and after
changing the Kummer generator of the torsor $f_K$, we can assume that
this torsor is given by an equation $X^p=1+\pi ^{n'} v$, where 
$v\in A'$ does not
belong to the ideal $\pi A'$, and $v$ is not a $p$-power mod $\pi$. 
Also one can check, after localisation and completion at the ideal 
$(\pi, S')$ as above (namely by looking what happend above 
the boundary $\Cal X_1$), that
necessarily $n'\le pv_K(\lambda)$, and $n'=n p$ is divisible by $p$
(cf. for example [Gr-Ma] proof of 1.1, and [He] proof of 1.6, chap 5). Moreover
after changing the Kummer generator of the torsor $f_K$, one can assume that
$\bar v={t'}^m\bar v'$, where $m$ is an integer prime to $p$, and $v'$ is 
a unit. Let $T:=T'({v/T'}^m) ^{1/m}$. We have $\Cal X\simeq 
\Spf R[[S,T]]/(ST-\pi ^{tp})$ for a suitable $S\in A'$.  
The cover $f:\Cal Y\to \Cal X$ is thus generically 
given by an equation $X^p=1+\pi ^{np}T^m$, and 
we are in case d or e. Let $\Cal X_1$ be 
as above the boundary of $\Cal X$ corresponding to the ideal $(\pi, T)$. 
The torsor $f_1:\Cal Y_1\to \Cal X_1$ induced by $f$ above $X_1$ is 
given by an equation $X^p=1+\pi ^{p(n+t)}S^{-m}$, which imply 
that $t+n\le v_K(\lambda)$ since $f_1$ is not completely split. Moreover
we are in case $d$ if $t+n=v_K(\lambda)$ and in case $e$ if 
$n+t<v_K(\lambda)$.

\pop {.5}
\par
\noindent
{\bf \gr 3.3. Variation of the different.}\ \rm (Compare with [He], 5.2) 
The following lemma, which is a direct consequence of 3.2.5, 
describes how does the degree of the different vary from one boundary to 
another in a cover $f:\Cal Y\to \Cal X$ as in 3.2.5.

\pop {.5}
\par
\noindent
{\bf \gr 3.3.1. Proposition.}\ \rm {\sl Let $\Cal X$ be the formal germ of a
semi-stable $R$-curve at an ordinary double point $x$. Let $f:\Cal Y\to \Cal X$
be a Galois cover of degree $p$, with $\Cal {Y}_k$ reduced and local, and 
with $f_K:\Cal Y_K\to \Cal X_K$ {\bf \'etale}. 
Let $y$ be the closed point of $\Cal Y$. Assume that $\bold {g_y=0}$,
which imply necessarily that the thikness $e=pt$ of the double 
point $x$ is divisible by $p$. For each integer $0<t'<t$, 
let $\Cal X_{t'}\to \Cal X$ be the blow-up of $\Cal X$ at the ideal 
$(\pi ^{pt'},T)$. The special fibre of $\Cal X_{t'}$ consists of a 
projective line $P_{t'}$ which meets two germs of double points $x$ and 
$x'$. Let $\eta$ be the generic point of  $P_{t'}$, and let $v_{\eta}$ 
be the corresponding discrete valuation of the function field of $\Cal X$.
Let $f_{t'}:\Cal Y_{t'}\to \Cal X_{t'}$ be the pull back of $f$, 
which is a Galois cover of degree $p$, and let $\delta (t')$ be the degree of 
the different induced by this cover above  $v_{\eta}$ (cf. [Sa] 3.1).
Also let ${\Cal X_i}$, for $i=1,2$, be 
the boundaries of $\Cal X$. Let $f_i:\Cal Y_i\to \Cal X_i$ be the 
torsors induced by $f$ above $\Cal X_i$, let $(G_{k,i},m_i,h_i)$ be their 
degeneration type, and let $\delta _i$ be the corresponding degree of the
different (cf. [Sa] 3.1). Say $\delta _1=\delta (0)$, 
$\delta _2=\delta (t)$, and $\delta (0)\le \delta (t)$. We have
$m:=-m_1=m_2$ say is positif. Then the 
following holds:
\par
1)\ If $\delta (t')<v_K(p)$, for every $0\le t'\le t$. Then for 
$0\le t_1\le t_2\le t$ we have $\delta (t_2)=\delta (t_1)+m(p-1)(t_2-t_1)$,
and $\delta (t')$ is an increasing function of $t'$.
\par
2)\ If there exists $0\le t'\le t$ such that $\delta (t')=v_K(p)$. Then there 
exists $0\le t_1 \le t_2\le t$ such that $\delta (t')=v_K(p)$ is constant
for $t_1\le t'\le t_2$, $\delta (t')$ is increasing as $t'$ 
increases from $0$ to $t_1$, and $\delta (t')$ decreases as $t'$
increases from $t_2$ to $t$.

\pop {.5}
\par
\noindent
{\bf \gr 3.3.2. Remark.}\rm \ In [Gr-Ma] and [He] where studied order 
$p$-automorphisms of open $p$-adic discs and $p$-adic annuli. In their aproach
one writes such an automorphism as a formal series and one deduce 
some results, e.g. the only if part of 3.1.2 and 3.2.3,  using 
the Weirstrass preparation theorem. The aproach adobted hier, and which 
consists on directly computing the 
vanishing cycles first, I believe provides another way to study such 
automorphisms. Namely these are those covers above formal fibres of 
semi-stable $R$-curves with genus $0$, and one can easily writes down  
Kummer equations which lead to such covers as in 3.1.3, 3.2.4, and 3.2.5.

\pop {2}
\par
\noindent
{\bf \gr References}
\rm

\pop {.5}
\par
\noindent
[Bo-Lu] S. Bosch and W. L\"utkebohmert, {\sl Formal and rigid geometry}, 
Math. Ann, 295, 291-317, (1993).

\pop {.5}
\par
\noindent
[Bo-Lu-1] S. Bosch and W. L\"utkebohmet, {\sl Stable reduction and 
uniformization of abelian varieties I}, Math. Ann. 270, 349-379, (1985). 

\pop {.5}
\par
\noindent
[Bo] N. Bourbaki, Alg\`ebre commutative, chap. 9, Masson, 1983.

\pop {.5}
\par
\noindent
[Gr-Ma] B. Green M. and M. Matignon, {\sl Order $p$-automorphisms of the 
open disc of a $p$-adic field}. J. Amer. Soc. 12 (1), 269-303, (1999). 

\pop {.5}
\par
\noindent
[Gr] A. Grothendieck, S\'eminaire de g\'eom\'etrie alg\'ebrique SGA-1, 
Lecture Notes 224, Springer Verlag, (1971).

\pop {.5}
\par
\noindent
[He] Y. Henrio, {\sl Arbres de Hurwitz et automorphismes d'ordre $p$ des 
disques et couronnes $p$-adiques formels}. Th\`ese de doctorat, 
Bordeaux France (1999).

\pop {.5}
\par
\noindent
[Ka] K. Kato, {\sl Vanishing cycles, ramification of valuations, and class 
field theory}, Duke Math. J. 55, 629-659, (1987).

\pop {.5}
\par
\noindent
[Ma-Yo] M. Matignon and T. Youssefi, {\sl Prolongement de morphismes de
fibres formelles}, Manuscript, (1992).

\pop {.5}
\par
\noindent
[Pr] R. J. Pries, {\sl Construction of covers via formal and rigid geometry}. 
In Courbes semi-stables et groupe fondamental en 
g\'eom\'etrie alg\'ebrique,
(J-B. Bost, F. Loeser, M. Raynaud, \'ed), Progress in Math. Vol 187, (2000).

\pop {.5}
\par
\noindent
[Ra] M. Raynaud, {\sl Rev\^etements de la droite affine en 
caract\'eristique $p>0$ et conjecture d' Abhyankar}, Invent. Math. 116, 
425-462 (1994).

\pop {.5}
\par
\noindent
[Ra-1] M. Raynaud, {\sl $p$-Groupes et reduction semi-stable des courbes},
The Grothendieck Festschrift, vol. 3, 179-197, Birkh\"auser, (1990).

\pop {.5}
\par
\noindent
[Ra-2] M. Raynaud, {\sl Sp\'ecialisation des rev\^etements en caract\'eristique
$p>0$}, Ann. Scient. Ec. Norm. Sup. 22, 345-375, (1989).

\pop {.5}
\par
\noindent
[Sa] M. Sa\"\i di, {\sl Torsors under finite and flat group schemes 
of rank $p$ with Galois action}, preprint.

\pop {.5}
\par
\noindent
[Sa-1] M. Sa\"\i di, {\sl Galois covers of degree $p$: semi-stable 
reduction and Galois action}, preprint.

\pop {.5}
\par
\noindent
[Sai] T. Saito, {\sl Vanishing cycles and differentials of curves 
over a discrete valuation ring}, Advanced Studies in Pure Mathematics 12, 
(1987).

\pop {2}
Mohamed Sa\"\i di
\pop {1}
\par
Departement of Mathematics
\par
University of Durham
\par
Science Laboratories
\par
South Road
\par
Durham, DH1 3LE, UK
\par
saidi\@durham.ac.uk

\enddocument
\end